\documentclass{elsart}
\usepackage{amssymb}
\usepackage{bm}
\usepackage{algorithmic}
\usepackage{natbib}

\newcommand{\krnl}{{\rm{Ker}}}
\newcommand{\ka}{{{k}}}
\newcommand{\expo}{{\rm{exp\: }}}

\newlength{\resto}
\newcommand{\lres}[2]%
 {%
 \settoheight{\resto}{$#1$}%
 \addtolength{\resto}{-1pt}%
 \raisebox{\resto}{\scriptsize $#2$}\overline{#1}%
 }
\newcommand{\pmb}[1]{{\boldsymbol{#1}}}
\newcommand{\Expo}{{\rm{Exp\: }}}
\newcommand{\INIT}{{\bf{INITIALIZATION:\ }}}
\newcommand{\env}{{{R^{\rm{env}}}}}
\newcommand{\level}{{\rm{level\: }}}
\newcommand{\estable}{{{{\mathbb{N}}^{n,(s)}}}}

\begin{document}
\begin{frontmatter}
\title{Gr\"obner bases and syzygies on bimodules}
\author{Manuel Garc\'{\i}a Rom\'an} \and
\author{Socorro Garc\'{\i}a Rom\'an}
\address{Departamento de Matem\'atica Fundamental,
Universidad de La Laguna, \\ E-38271 La Laguna, Spain}
\thanks{Research supported by PI2002/075, Gobierno de Canarias.}

\begin{abstract}
A new more efficient method for the computation of two-sided Gr\"obner bases of
ideals and bimodules shifting the problem to the enveloping algebra is proposed.
Arising from the ideas this method involves, we introduce the notion of two-sided
syzygy, which reveals to be useful in the computation of the intersection of bimodules.
Further applications are left for a sequel.

\end{abstract}
\end{frontmatter}


\begin{section}{Introduction}

Though first developed in the ring of polynomials, the methods based on Gr\"obner
bases also work in some noncommutative rings, e.g. the Weyl algebras or, more generally,
the so-called Poincar\'e-Birkhoff-Witt rings (PBW, for short), including some classical
quantum groups. After first results in
the Weyl algebra \citep{galligo} and in tensor algebras of finite-dimensional Lie algebras
\citep{apel}, Kandri-Rodi and Weispfenning were the first to introduce Gr\"obner bases
in the more general class of algebras where the degree of a skew-commutator $p_{ij}=x_jx_i-c_{ij}x_ix_j$ is bounded
by the degree of the product of generators $x_ix_j$, for $1\le i<j\le n$
\citep[see][]{KW}. The theory has recently been surveyed in \citet{BCGL98}.
Algorithms to compute the Gelfand-Kirillov dimension, to check whether a two-sided
ideal is prime or not and to compute the projective dimension of a module have also been developed \citep[see][]{BCGL96,Lob,BCGL99,Gag}.

On these generalizations, authors were mainly interested in one-sided
ideals and modules, whereas methods for the two-sided counterparts are merely patches
in order to cope with the two-sided input data~\citep[cf.][]{Pes,BGV}

In this note we show that those mends are not necessary, due to the very well known
fact that two-sided ideals and bimodules may be seen as left modules on
the enveloping algebra. First, we show that the enveloping algebra of a PBW algebra
is another PBW algebra. Second, we find a method to shift the data back and forth
through the morphism
 \[ \mathfrak{m}^s : (R^{\rm{env}})^s \longrightarrow R^s \ ; \ \
 (f_i\otimes g_i)_{i=1}^{s} \longmapsto (f_ig_i)_{i=1}^{s}\]
in order to carry out the computations on the enveloping algebra using one-sided techniques.

This philosophy allows, for example, to compute Gr\"obner bases for
bimodules with only one call to the left Buchberger algorithm, instead of the a
priori unknown number of calls typical of the aforementioned methods.

The techniques we use have led us to study the syzygy bimodule, which is the
two-sided counterpart of the left syzygy module. Amongst its applications,
we show that it can be used in the
computation of intersections of bimodules when one starts, as usual, from two-sided
input data.


Throughout this paper we will use the following notation.
We will denote by
$\epsilon_i$ the element
$(0,...,\stackrel{i}{1},...,0) \in \mathbb{N}^n$.
The symbol $x^{\alpha}$ will denote the monomial $x_1^{\alpha_1}\cdots x_n^{\alpha_n}$ in the free algebra $\ka \langle x_1, \ldots, x_n \rangle$ or in any of its epimorphic images, $\ka$ being a field.
$R$ will be a PBW $\ka$-algebra, $R^{\rm{op}}$ its
opposite algebra and $\env$ its en\-ve\-lo\-ping
algebra $R \otimes_{\ka} R^{\rm{op}}$.
Finally, for any subset $F$ of the free left $R$-module $R^s$, we will denote by
$_R\langle F \rangle $, resp.\ $_R\langle F \rangle_R$, the left $R$-module, resp.\
the $R$-bimodule, generated by $F$.

The computations of the examples shown in this paper were done using a library
of procedures built by the
authors using the package of symbolic computation Maple 6.
The computation times correspond to a Pentium III 700 Mhz personal computer with 256 Mb RAM.

\end{section}


\begin{section}{The enveloping algebra of a PBW algebra}

A subset $Q= \{ x_jx_i- q_{ji}x_ix_j - p_{ji};\ 1 \leq i < j \leq n\}$ of
the free algebra $k\langle x_1,\ldots,x_n \rangle$
is a set of
quantum relations bounded
by the admissible order ``$\preceq$'' on $\mathbb{N}^n$ if
$q_{ji} \in k^{\ast}$ and $p_{ji}$ is a finite $\ka$-linear combination of
standard monomials $x^{\alpha}$ ($\alpha \in \mathbb{N}^n$) such that
$\expo(p_{ji})\prec \epsilon_i +\epsilon_j$, for all $i<j$, where $\expo(f)$
denotes the exponent of the leading term of the element $f$.
As may be found in the literature, \citep[see e.g.][]{KW, BCGL98, BGV}, a
{\em{Poincar\'e-Birkhoff-Witt algebra}\/}
(PBW algebra, for short) is a $\ka$-algebra $R$
where the set of standard monomials $\{x^{\alpha}; \: \alpha \in \mathbb{N}^n\}$
is a $\ka$-basis and such that
there exists a set of quantum relations $Q$ bounded by an ad\-mi\-ssi\-ble
order ``$\preceq$'' satisfying
\[
R=\frac{k\langle x_1, \ldots, x_n\rangle}{_{k\langle x_1, \ldots, x_n\rangle}\langle
\: Q\: \rangle_{k\langle x_1, \ldots, x_n\rangle} }.
\]
This algebra is usually denoted by $k\{x_1, \ldots, x_n;\: Q, \preceq\}$.

Amongst the examples of PBW algebras, we find the commutative
polynomial ring $k[x_1,...,x_n]$, some iterated Ore extensions or
a pretty large class of quantum groups just as the
multiparameter $n$-dimensional quantum space $\mathcal{O}_q (\mathbb{A}^n)$,
the bialgebra of quantum matrices $M_q(2)$,
the Weyl algebra $A_n(k)$,
the enveloping algebra of traceless matrices $U(\mathfrak{sl}(2))$, etc.

The tensor product of PBW algebras is a new PBW algebra:

\begin{prop}\label{productotensorialalgebrasPBW}
If $R=k\{x_1, \ldots,x_m;\: Q_{R}, \preceq_{R}\}$ and
$S=k\{y_1,\ldots,y_n;\: Q_{S}, \preceq_{S}\}$ are PBW algebras
with quantum relations
\[
Q_{R}= \{ x_jx_i - q_{ji}x_ix_j - p_{ji};\: 1 \leq i < j \leq m\},
\]
\[
Q_{S}= \{ y_jy_i - q'_{ji}y_iy_j - p'_{ji};\: 1 \leq i < j \leq n \},
\]
then $R \otimes_k S$ is the PBW algebra denoted by
\[
k\{x_1\otimes 1,\ldots,x_m\otimes 1,
1\otimes y_1,\ldots,1\otimes y_n; Q, \preceq\}, {\rm{\ \ where}}
\]
\[
Q=\left\{\begin{array}{cc} (x_j \otimes 1)(x_i \otimes 1)-
q_{ji}(x_i \otimes 1)(x_j \otimes 1)- p_{ji} \otimes 1;
& 1 \leq i < j \leq m \\
(1 \otimes y_j)(x_i \otimes 1) - (x_i \otimes 1)(1 \otimes y_j);
& 1 \leq i \leq m,\; 1 \leq j \leq n \\
(1 \otimes y_j)(1 \otimes y_i)-
q'_{ji}(1 \otimes y_i)(1 \otimes y_j) - 1 \otimes p'_{ji};
& 1 \leq i < j \leq n.
\end{array} \right\}
\]
and ``$\preceq$" is one amongst the elimination orders \citep[as defined, e.g.,
in][]{CLO}
arising from ``$\preceq_{R}$"
and ``$\preceq_{S}$".
\end{prop}

As a first example, note that $A_{n+m}(\ka)$ is the PBW algebra
$A_n(\ka) \otimes A_m(\ka)$ constructed in the proposition.
Another example of this construction
is the enveloping algebra $\env$
of $R= k\{x_1, \ldots, x_n;\: Q, \preceq\}$. Before we describe it, let us define
the composition orders. For any $\alpha =(\alpha_1,...,\alpha_1)\in \mathbb{N}^n$,
denote by $\alpha^{\rm{op}}$ the $n$-tuple $(\alpha_n,...,\alpha_1)$.

\begin{defn}\label{compositionorder}
Let ``$\preceq$" be an order on $\mathbb{N}^n$.
The {\em{up-component}\/} composition order in $\mathbb{N}^{2n}$, denoted
 ``$\preceq^{c}$", is defined by
\[
(\alpha,\beta) \prec^c (\gamma, \delta) \Leftrightarrow
\left\{ \begin{array}{l}
\alpha + \beta^{\rm{op}}  \prec  \gamma + \delta^{\rm{op}} ,\ 
 {\rm{or}}  \\
\alpha + \beta^{\rm{op}}  =  \gamma + \delta^{\rm{op}}
{\rm{\ and\  }} \beta^{\rm{op}}  \prec  \delta^{\rm{op}}
\end{array} \right.
\]
The {\em{down-compontent}\/} composition order ``$\preceq_{c}$" is
defined by
\[
(\alpha,\beta) \prec_c (\gamma, \delta) \Leftrightarrow
\left\{ \begin{array}{l}
\alpha + \beta^{\rm{op}}  \prec  \gamma + \delta^{\rm{op}} ,\ 
 {\rm{or}}  \\
\alpha + \beta^{\rm{op}}  =  \gamma + \delta^{\rm{op}}
{\rm{\ and\  }} \alpha  \prec  \gamma
\end{array} \right.
\]
\end{defn}

If ``$\preceq$'' is an admissible order on $\mathbb{N}^n$, then both composition
orders ``$\preceq^c$'' and ``$\preceq_c$''are admissible orders on $\mathbb{N}^{2n}$.

Note that the opposite algebra $R^{\rm{op}}$ is the PBW
algebra ${k\{x_n,\ldots,x_1;\: Q^{\rm{op}}, \preceq^{op}\}}$, where the elements of
$Q^{\rm{op}}$ are those of $Q$ written oppositely and ``$\preceq^{op}$" is the order in
$\mathbb{N}^n$ given by $\alpha \preceq^{\rm{op}} \beta \iff \alpha^{\rm{op}} \preceq \beta^{\rm{op}}$.


\begin{prop}\label{productotensorialRRopuestoPBW}
If $R= k\{x_1,\ldots,x_n;\: Q, \preceq\}$ is a PBW algebra
with quantum relations $Q=\{x_jx_i-q_{ji}x_ix_j-p_{ji};\: 1 \leq i < j \leq n\}$,
then $\env$ is the PBW algebra
$k\{x_1\otimes 1,\ldots,x_n\otimes 1,
1\otimes x_n,\ldots,1\otimes x_1;\: Q^{\ast}, \preccurlyeq \},
$,
where
\[
Q^{\ast}=\left\{\begin{array}{cc} (x_j \otimes 1)(x_i \otimes 1)-
q_{ji}(x_i \otimes 1)(x_j \otimes 1)- p_{ji} \otimes 1;
\ & 1 \leq i < j \leq n \\
(1 \otimes x_j)(x_i \otimes 1) -(x_i \otimes 1)(1 \otimes x_j);
\  & 1 \leq i,j \leq n \\
(1 \otimes x_i)(1 \otimes x_j)-q_{ji}(1 \otimes x_j)(1 \otimes x_i)
- 1 \otimes p_{ji};\  & 1 \leq i < j \leq n.
\end{array} \right\}
\]
and ``$\preccurlyeq$" is, either any of the elimination orders ``$\preceq^{\ast}$"
or ``$\preceq_{\ast}$" in $\mathbb{N}^{2n}$
co\-rres\-pon\-ding to ``$\preceq$"
and ``$\preceq^{\rm{op}}$",
or any of the composition orders ``$\preceq^{c}$" or ``$\preceq_{c}$" on $\mathbb{N}^{2n}$
corresponding to ``$\preceq$".
\end{prop}

In what follows,
we will work on the free $R$-bimodule $R^s$, being $s$ a positive integer and
$R$ a PBW algebra, and we will use the $R$-module basis $\{{\pmb{e_i}}\}_{i=1}^s$
consisting of ${\pmb{e_i}}=(0,...\stackrel{i}{1}...,0) \in R^s$ for all $1 \leq i \leq s$.

The notion and some applications of left Gr\"obner bases in PBW algebras and left
modules may be found, e.g., in \citet{BGV}.
For convenience, just recall that
\label{baseGrobnerBilatera}
if $M \subset R^s$ is an $R$-bimodule, then $G=\{{\pmb{g_1}},...,\: {\pmb{g_r}}\} \subset M \setminus \{0\}$ is a two-sided Gr\"obner basis if one of the following
equivalent statements holds
\begin{enumerate}
\item $M=\: _R\langle G \rangle_R$ and
$\Expo(M)=\bigcup_{k=1}^r (\mathbb{N}^n+\expo(\pmb{g_k}))$;

\item $G$ is a left Gr\"obner basis and $M=\: _R\langle G \rangle_R=\:
_R\langle G \rangle$;

\item $G$ is a left Gr\"obner basis, $M=\: _R\langle G \rangle_R$ and
${\pmb{g_k}}x_i \in\: _R\langle G \rangle$, for all $k \in \{1,...,r\}$ and $i \in \{1,...,n\}$.
\end{enumerate}
A set $G \subset R^s$ is said to be a {\em{two-sided Gr\"obner basis}\/} if
so it is for the $R$-bimodule $_R\langle G \rangle_R$.
\end{section}

\begin{section}{Computing two-sided Gr\"obner Bases}

This section is devoted to the methods for computation of two-sided Gr\"obner
bases for $R$-subbimodules of $R^s$.
We denote by ${\pmb{f}} \otimes {\pmb{g}}$ the element $(f_1 \otimes g_1, ...,f_s \otimes g_s) \in (\env)^s$, where
${\pmb{f}}=(f_1,...,f_s), {\pmb{g}}=(g_1,...,g_s) \in R^s$.

As a consequence of the third caracterization
of two-sided Gr\"obner bases above,
some authors have proposed an algorithm to compute them \citep[see][]{BGV}.
Alternatively, we propose a new algorithm which improves that one,
since it calls only once the left Buchberger Algorithm, although it uses
more variables and input elements. The philosophy is to transform
the problem into computing a left Gr\"obner basis in the free module
$(\env)^s$. This may be done since, just as we saw in the previous
section, $\env$ has a PBW structure.

It is known that $R$-bimodules are exactly left $\env$-modules.
%
Notice that, in particular, the free module $R^s$ is a left $\env$-module
with the action
$(r \otimes r')\: {\pmb{f}}=(rf_1r',...,rf_sr')$, and
$(\env)^s$ possesses an $R$-bimodule structure
whose multiplicacions are given by
$r({\pmb{f}} \otimes {\pmb{g}})r'=(rf_1 \otimes g_1r',...,rf_s \otimes g_sr')$, where
$r, r' \in R$ and ${\pmb{f}}=(f_1,...,f_s),\ {\pmb{g}}=(g_1,...,g_s) \in R^s$.

Likewise, the map $\mathfrak{m}^s=\mathfrak{m} \times \stackrel{s}{\cdots}
\times \mathfrak{m}: (\env)^s \rightarrow R^s$, where $\mathfrak{m}(r \otimes r')=rr'$, for $r, r' \in R$, is
an epimorphism of left $\env$-modules. Thus there exists a bijection \\
\[
\begin{array}{rcl}
\{N \subseteq (\env)^s ;\: \krnl ({\mathfrak{m}^s}) \subseteq N \in \env{\rm{-Mod}}\} & \longrightarrow &
\{M \subseteq R^s ;\: M \in R{\rm{-Bimod}}\} \\
N & \rightarrow & M_N:=\mathfrak{m}^s(N), \\
N_M:=(\mathfrak{m}^s)^{-1}(M) & \leftarrow & M
\end{array}
\]

Using this bijection, for each $R$-bimodule $M \subset R^s$, we have a left
$\env$-module $N_M \subset (\env)^s$. Moreover, from a finite generator system for $M$
it is possible to obtain one for $N_M$, just as
the following results show.
\begin{lem}
Let $R$ be a $\ka$-algebra.
\begin{enumerate}
\item
If $M=\: _R \langle{\pmb{f_1}},...,\: {\pmb{f_t}} \rangle_R \in R^s$, then
$N_M=\: _{\env} \langle {\pmb{f_1}} \otimes {\pmb{1}},...,\: {\pmb{f_t}} \otimes {\pmb{1}} \rangle + \krnl(\mathfrak{m}^s)$;

\item
$\krnl(\mathfrak{m}^s)=\: _{\env} \langle {\pmb{f}} \otimes\: {\pmb{1}} - {\pmb{1}} \otimes {\pmb{f}}
\: ;\: {\pmb{f}} \in R^s \rangle$;

\item
If $R=k\{x_1,...,x_n;\ Q;\ \preceq\}$ is a PBW algebra, then
$_{\env}\: \langle {\pmb{f}} \otimes {\pmb{1}} - {\pmb{1}} \otimes {\pmb{f}} \: ;\: {\pmb{f}} \in R^s \rangle
=\: _{\env}\langle {\pmb{x^{(\epsilon_j,k)}}} \otimes {\pmb{1}}-{\pmb{1}} \otimes {\pmb{x^{(\epsilon_j,k)}}}
\: ;\: 1\leq j \leq n, 1 \leq k \leq s \rangle$.
\end{enumerate}
\end{lem}

\begin{cor}\label{corolariogeneradoresMyN}
If $R=k\{x_1,...,x_n;\ Q;\ \preceq\}$ is a PBW algebra and $M=\: _R \langle {\pmb{f_1}},...,\: {\pmb{f_t}} \rangle_R$ $\subseteq R^s$ then

\[ N_M=\: _{\env} \langle\ \{{\pmb{f_i}} \otimes {\pmb{1}}\}_{i=1}^t, \{{\pmb{x^{(\epsilon_j,k)}}} \otimes {\pmb{1}}-
{\pmb{1}} \otimes {\pmb{x^{(\epsilon_j,k)}}}\}_{1\leq j \leq n,\ 1 \leq k \leq s}\ \rangle. \]
\end{cor}

Given an admissible order in $\mathbb{N}^n$, from here on we will call TOP
(term over position), resp.\ POT (position over term) the orders in $\estable$
given by
\[
(\alpha,i) \prec (\beta,j) \iff
\left\{ \begin{array}{l}
\alpha \prec \beta, \ {\rm{or}}  \\
\alpha = \beta {\rm{\ and\  }} i>j
\end{array} \right.
\ \ {\rm{resp.}} \ \ \
\left\{ \begin{array}{l}
i>j , \ {\rm{or}}  \\
i = j {\rm{\ and\  }} \alpha \prec \beta
\end{array} \right.
\]

\begin{lem}
Let $R=k\{x_1,...,x_n;\ Q;\ \preceq\}$ be a PBW algebra and consider the order TOP (or POT) on both
$R^s$ and $(\env)^s$.
\begin{itemize}
\item Taking $\preceq^{\ast}$ or $\preceq^c$ on $\env$, if ${\pmb{h}} \in (\env)^s$ is such that
$\expo_{(\env)^s}({\pmb{h}})=((\alpha, 0),i) \in \mathbb{N}^{2n,(s)}$, then ${\pmb{h}} \not\in \krnl(\mathfrak{m}^s)$ and
$\expo_{R^s}(\mathfrak{m}^s({\pmb{h}}))=(\alpha,i)$;

\item Taking $\preceq_{\ast}$ or $\preceq_c$ on $\env$, if ${\pmb{h}} \in (\env)^s$ is such that
$\expo_{(\env)^s}({\pmb{h}})=((0,\alpha),i) \in \mathbb{N}^{2n,(s)}$, then ${\pmb{h}} \not\in \krnl(\mathfrak{m}^s)$ and
$\expo_{R^s}(\mathfrak{m}^s({\pmb{h}}))=(\alpha^{\rm{op}},i)$.
\end{itemize}
\end{lem}

Using these results we have:
\begin{thm}\label{algoritmoBGBilateras}
Let $R=k\{x_1,...,x_n; Q, \preceq\}$ be a PBW algebra, $M \subset R^s$ be an $R$-bimodule
and consider in $\env$ the PBW structure given in the proposition
\ref{productotensorialRRopuestoPBW}
(where the order is one of  $\preceq^{\ast}$, $\preceq^c$, $\preceq_{\ast}$ or $\preceq_c$).

If $G$ is a left Gr\"obner basis for $N_M=(\mathfrak{m}^s)^{-1}(M)$ with
TOP (resp.\ POT), then the set $\mathfrak{m}^s(G) \setminus \{0\}$ is a
two-sided Gr\"obner basis for $M$ with TOP (resp.\ POT).
\end{thm}
The theorem provides a method of construction of two-sided Gr\"obner bases
for bimodules. For convenience we write explicitly under the name of
Algorithm~\ref{algoritmobaseGrobnerbilatera}.

\begin{alg}
Two-sided Gr\"obner bases 
\begin{algorithmic}\label{algoritmobaseGrobnerbilatera}
\smallskip
\REQUIRE $F=\{{\pmb{f_1}},...,\: {\pmb{f_t}}\}\subseteq R^s \setminus \{0\}$;
\ENSURE $G=\{{\pmb{g_1}},...,\: {\pmb{g_{t'}}}\}$, a two-sided Gr\"obner basis for
$\: _{R}\langle F \rangle_{R}$ such that $F \subseteq G$;

\medskip
\INIT $B:= \{{\pmb{f_i}} \otimes {\pmb{1}}\}_{i=1}^t \cup \{{\pmb{x^{(\epsilon_j,k)}}} \otimes {\pmb{1}}-
{\pmb{1}} \otimes {\pmb{x^{(\epsilon_j,k)}}}\}_{1\leq j \leq n,\ 1 \leq k \leq s}$;
\vskip0.1cm
\STATE  Using the Left Buchberger Algorithm, compute a left Gr\"obner basis $G'$
in the PBW algebra $(\env)^s$ for the input data $B$;
\medskip
\STATE If $G'=\{{\pmb{g'_1}},...,\: {\pmb{g'_r}}\}$ with ${\pmb{g'_i}}= (\sum_{j \in \mathfrak{I}_i}
p^1_{ij} \otimes q^1_{ij},...,
\sum_{j \in \mathfrak{I}_i} p^s_{ij} \otimes q^s_{ij})$, take

${\pmb{g_i}}:= (\sum_{j \in \mathfrak{I}_i} p^1_{ij} q^1_{ij},..., \sum_{j \in \mathfrak{I}_i} p^s_{ij} q^s_{ij})$;
\medskip
\STATE $G:= \emptyset$;
\medskip
\FORALL{$i=1$ to $r$}
  \IF{${\pmb{g_i}} \neq 0$}
    \STATE $G:= G \cup \{{\pmb{g_i}}\}$;
  \ENDIF
\ENDFOR
\medskip
\end{algorithmic}
\end{alg}

The advantage offered by this algorithm is that only one call to the left Buchberger Algorithm
is done, whereas
the one shown in \citep{BGV} makes an, a priori, unknown number of calls.
\begin{exmp}\label{planocuantico}
Let $R$ be the quantum plane, i.e.,
$R={\mathbb{C}\{x,y;\ \{yx-qxy\}, \preceq_{(1,3)}\}}$
(take $q=i\sqrt{2}$) where $\preceq_{(1,3)}$ is the $(1,3)$-weighted
lexicographical order.
Let
$F=\{(2x,x^2y,xy^2+y^2),(xy,0,-x^2y^2),(x^2,2,0)\} \subset R^3$
and consider the order TOP in $R^3$.

The old algorithm \citep[see][]{BGV} takes 66.371 seconds to compute the two-sided
Gr\"obner basis, consisting of 17 elements:
\begin{eqnarray*}
\lefteqn{G =  \{\: (2x,x^2y,xy^2+y^2),\: (xy,0,-x^2y^2),\: (x^2,2,0),\: (xy-2x,x^3y-x^2y-4,-y^2),} \\
& & (xy,x^4y+y-4x,0),\: (0,(i\sqrt{2}+2)x^3y-12,0),\: (0,(i\sqrt{2}+2)y,0), \\
& & ((-i\sqrt{2}-2)xy,(-i\sqrt{2}-2)x^3y+(-4+4i\sqrt{2})x+12,0),\: (0,(-4+4i\sqrt{2})x^2,0), \\
& & ((i\sqrt{2}-1)xy^2+(2-2i\sqrt{2})xy,(2-2i\sqrt{2})x^3y^2+3x^2y^2+(-4-2i\sqrt{2})y,0), \\
& &  ((1-i\sqrt{2})xy^2+(-2+2i\sqrt{2})xy, (1+2i\sqrt{2})x^3y^2-3x^2y^2,0), \\
& &  ((1+i\sqrt{2}/2)xy^2,0,0),\: (0,3y,0),\:((1-i\sqrt{2}/4)xy^2,3/4y^2+(-4+i\sqrt{2})xy,0),\\
& & (xy,-2i\sqrt{2}x,0),\:(0,(-4+2i\sqrt{2})x,0), \: (0,-4+2i\sqrt{2},0)\: \}.
\end{eqnarray*}
The algorithm \ref{algoritmobaseGrobnerbilatera} takes 46.972 seconds to
compute a two-sided Gr\"obner basis with 12 elements
(considering the order $\preceq^{c}$ (see remark \ref{compositionorder}) in $\env$):
\begin{eqnarray*}
\lefteqn{G = \{\: (2x,x^2y,xy^2+y^2),\: (xy,0,-x^2y^2),\: (x^2,2),\: (xy-2x,x^3y-x^2y-4,-y^2),} \\
& & (0,(i\sqrt{2}+2)x^3y-12,0),\: ((i\sqrt{2}-1)xy^2+(2-2i\sqrt{2})xy,3x^2y^2+(2-2i\sqrt{2})y, \\
& & 0),\: (2xy,-i\sqrt{2}y-4i\sqrt{2}x,0),\: (0,(-i\sqrt(2)/4-1/2)y^2+(2-2i\sqrt{2})xy,0),\\
& & (0,3y,0),\: (0,-2i\sqrt{2}x^2,0),\: (0,-6i\sqrt{2}x,0),\: (0,-16i\sqrt{2}-16,0)\: \}.
\end{eqnarray*}
\end{exmp}

\begin{exmp}
Let $R$ be the PBW algebra
$M_q(2)=\mathbb{C} \{ x, y, z, t;  Q, \preceq_{\rm{glex}} \}$
 of quantum matrices (with $q=i\sqrt{2}$), where
$ Q=\{
yx-qxy, 
ty-qyt, 
zx-qxz, 
tz-qzt, 
zy-yz, 
tx-xt-(q^{-1}-q)yz
\}
$
and consider the order ``POT" in $R^2$.
Let $F=\{(-xzt+1,2y^3),(x^2t,y^2)\} \subset R^2$.

The old algorithm takes 216.433 seconds to compute the two-sided Gr\"obner basis,
consisting of 26 elements:
\begin{eqnarray*}
\lefteqn{G = \{\: (-xzt+1,2y^3),\: (x^2t,y^2),\: (-x,-2xy^3+1/2y^2z),\: (-3/2yz^2-1,4i\sqrt{2}xy^3zt } \\
& & -3y^4z^2-y^2z^2t-2y^3),\: (0,4i\sqrt{2}x^2y^3t-i\sqrt{2}xy^2zt+3/8i\sqrt{2}y^3z^2+y^2),\\
& & (0,3i\sqrt{2}xy^4z^2+3/8y^3z^3-(4i\sqrt{2}+2)xy^3+(1/2-1/2i\sqrt{2})y^2z),\\
& & (0,-3i\sqrt{2}x^2y^4z-3/8i\sqrt{2}xy^3z^2-3xy^2),\: (0,(4i\sqrt{2}+2)x^2y^3\\
& & -(i\sqrt{2}+1/2)xy^2z),\: (0,(12-6i\sqrt{2}xy^4z^2-3/4(i\sqrt{2}+1)y^3z^3\\
& & +(2+4i\sqrt{2})xy^3+3/2y^2z),\: (0,-12x^2y^4z-3/2xy^3z^2+(-2+2i\sqrt{2})xy^2),\\
& & (0,(1+1/2i\sqrt{2})y^3),\:(0,(1/2+1/4i\sqrt{2})y^3z),\: (0,(1-i\sqrt{2})y^3z^2t),\: (0,\\
& & (i\sqrt{2}+2)xy^3zt+(3/8i\sqrt{2}-3/8)y^4z^2+(1+1/2i\sqrt{2})y^3),\: (0,(1+1/2i\sqrt{2})y^2z),\\
& & (0,(1/2+1/4i\sqrt{2})y^2z^2),\:(0, (1-i\sqrt{2})y^2z^3t),\: (0,(i\sqrt{2}+2)xy^2z^2t\\
& & +(3/8i\sqrt{2}-3/8)y^3z^3+(1+1/2i\sqrt{2})y^2z),\: ((1+1/4i\sqrt{2})t,-8xy^3zt^2\\
& & +3i\sqrt{2}y^4z^2t-i\sqrt{2}y^2z^2t^2+4y^3t),\: (3/2y,3y^4+3y^2t),\:(3/2i\sqrt{2}yz,\\
& & (-4i\sqrt{2}-2)xy^3t+3i\sqrt{2}y^4z+(1/2+i\sqrt{2})y^2zt),\:((-1+1/4i\sqrt{2})t,\\
& & (4i\sqrt{2}+8)xy^3zt^2+(-3-3i\sqrt{2})y^4z^2t+(i\sqrt{2}-1)y^2z^2t^2),\:(0,-3i\sqrt{2}xy^4zt\\
& & +3/8i\sqrt{2}y^3z^2t+(1+1/2i\sqrt{2})y^2t),\:(-1,0),\: (0,-i\sqrt{2})y^2t),\:(0,(7/9-4/9i\sqrt{2})y^2)\:\}.
\end{eqnarray*}
The algorithm \ref{algoritmobaseGrobnerbilatera} takes 59.909 seconds to compute
the two-sided Gr\"obner basis (considering the order $\preceq^{\ast}$ in $\env$), consisting of 11 elements:
\begin{eqnarray*}
\lefteqn{ G = \{\: (-xzt+1,2y^3),\:(x^2t,y^2),\:(-x,-2xy^3+1/2y^2z),\:(0,-3i\sqrt{2}xy^4z^2-3/8y^3z^3 } \\
& & +(4i\sqrt{2}+2)xy^3+(1/2i\sqrt{2}-1/2)y^2z),\:(-3/2yz^2t+(-1/2i\sqrt{2}-1)t,-6y^3t),\\
& & (0,-3i\sqrt{2}x^2y^4z-3/8i\sqrt{2}xy^3z^2-3xy^2),\: (0,(1+1/2i\sqrt{2})y^3),\\
& & (0,(1+1/2i\sqrt{2})y^2z),\:(3/2y,3y^2t),\:(-1,0),\:(0,y^2)\: \}.
\end{eqnarray*}
\end{exmp}
\end{section}
\begin{section}{Syzygy Bimodules}

In this section we study the notion of {\em{syzygy bimodule}\/} of a subset of
$R^s$, $R$ being a PBW algebra and $s \in \mathbb{N}^{\ast}$. This notion
can be viewed as the ana\-lo\-gous one of {\em{left syzygy module}\/}
for left modules, since it presents some similar properties.




There exists an algorithm \citep[see again][]{BGV} which computes a generator
system of the left syzygy module $Syz^l(F)$, provided a finite set of input data $F \subset R^s$ is given.
This algorithm is shown below 
and will be used within Algorithm \ref{algoritmoBimoduloSicigia}.
\begin{alg}
Left Syzygy Module
\begin{algorithmic}\label{algoritmoModuloSicigiaIzda}
\smallskip
\REQUIRE $F=\{{\pmb{f_1}},...,\: {\pmb{f_t}}\}\subseteq R^s \setminus \{0\}$;
\ENSURE $H$, a finite left generator system of $Syz^l(F)$;
\medskip

\INIT
Run the Left Buchberger Algorithm for the input data $F$ in order to compute:

- a left Gr\"obner basis $G= \{{\pmb{g_1}},...,\: {\pmb{g_r}}\} \subset R^s$ for $_R\langle F \rangle$,

- the elements $h_{ij}^k \in R$ such that $SP({\pmb{g_i, g_j}})=\sum_{k=1}^r h_{ij}^k {\pmb{g_k}}$
for all $1 \leq i < j \leq r$, and

- the matrix $Q \in M_{r \times t}(R)$ such that $({\pmb{g_1}},...,\: {\pmb{g_r}})=({\pmb{f_1}},...,\: {\pmb{f_t}})Q^t$;
\medskip
\FORALL{$1\leq i < j \leq r$}
\smallskip
  \IF{$\level(\expo({\pmb{g_i}})) = \level(\expo({\pmb{g_j}}))$}
\smallskip
    \STATE Compute $r_{ij}, r_{ji}$ such that $SP({\pmb{g_i}}, {\pmb{g_j}})=r_{ij}{\pmb{g_i}}-r_{ji}{\pmb{g_j}}$;
    \STATE Let ${\pmb{p_{ij}}}:=(0,...\stackrel{i}{r_{ij}}...,0)-(0,...\stackrel{j}{r_{ji}}...,0)-(h_{ij}^1,...,h_{ij}^r)$;
\smallskip
   \ENDIF
\ENDFOR
\medskip
\STATE Let $H:=\{p_{ij}Q\ {\rm{\ such\ that\ }} 1 \leq i < j \leq r, {\rm{\ and\ }}
\level(\expo({\pmb{g_i}})) = \level(\expo({\pmb{g_j}}))\}$.
\medskip
\end{algorithmic}
\end{alg}
\begin{defn}
Let ${\pmb{f_1}},...,{\pmb{f_t}} \in R^s$.
The {\em{syzygy bimodule}\/} of the matrix
\[
F=\left[\begin{array}{c}
{\pmb{f_1}} \\
\vdots \\
{\pmb{f_t}}
\end{array}\right]
\in M_{t \times s}(R),
\]
denoted by $Syz(F)$ or $Syz(\pmb{f_1},...,\: \pmb{f_t})$,
is the kernel of the homomorphism of
left $\env$-modules
$(\env)^t  \longrightarrow   R^s ; \ \
(h_1,...,h_t)  \longmapsto  \sum_{i=1}^t h_i{\pmb{f_i}}$.
\end{defn}

We can compute the syzygy bimodule of a matrix $F$ using again the techniques showed
in section 3, that is, we will {\em{move}\/}  the problem to the context of
the enveloping algebra in order to use the methods on the left side.

\begin{prop}
Let $M\subseteq R^s$ be a $R$-bimodule and $N=(\mathfrak{m}^s)^{-1}(M)$.
Let $\{{\pmb{h_1}},...,\: {\pmb{h_r}}\} \subset (\env)^{t+sn}$ be a generator system of
$Syz^l(\{\pmb{f_i} \otimes {\pmb{1}}\}_{i=1}^t, \{{\pmb{x^{(\epsilon_j,k)}}} \otimes {\pmb{1}}-
{\pmb{1}} \otimes {\pmb{x^{(\epsilon_j,k)}}}\}_{1\leq j \leq n,\ 1 \leq k \leq s})$ as a left $\env$-module.

Then $Syz({\pmb{f_1}},...,\: {\pmb{f_t}})=\: _R\langle \pi({\pmb{h_1}}),...,\pi({\pmb{h_r}}) \rangle_R$, where
$\pi$ is the projection homomorphism
$\pi: (\env)^t \times (\env)^{sn}  \longrightarrow  (\env)^t$.
\end{prop}

\begin{pf}
Notice that $\sum_{i=1}^t g_i({\pmb{f_i}} \otimes {\pmb{1}}) \in \krnl(\mathfrak{m}^s)$
for any $\pmb{g}=(g_1,...,g_t) \in Syz({\pmb{f_1}},...,\: {\pmb{f_t}})$.
So, there exists an element
${\pmb{g'}}=(g_{11}',...,g_{1s}',...,g_{n1}',...,g_{sn}')$ $\in (\env)^{sn}$ such that
$\sum_{i} g_i({\pmb{f_i}} \otimes {\pmb{1}})=\sum_{j,k}
g_{jk}'({\pmb{x^{(\epsilon_j,k)}}} \otimes {\pmb{1}}-{\pmb{1}} \otimes {\pmb{x^{(\epsilon_j,k)}}})$. Hence,
$({\pmb{g}},-{\pmb{g'}})$ is in $Syz^l(\{\pmb{f_i} \otimes {\pmb{1}}\}_{i=1}^t, \{{\pmb{x^{(\epsilon_j,k)}}} \otimes {\pmb{1}}-
{\pmb{1}} \otimes {\pmb{x^{(\epsilon_j,k)}}}\}_{1\leq j \leq n,\ 1 \leq k \leq s})$. From this point, the proof may easily be finished. \qed
\end{pf}
\begin{alg}
Syzygy Bimodule
\begin{algorithmic}\label{algoritmoBimoduloSicigia}
\smallskip
\REQUIRE $F=\{{\pmb{f_1}},...,\: {\pmb{f_t}}\}\subseteq R^s \setminus \{0\}$;
\ENSURE $H$, a finite generator system of $Syz(F)$ as an $R$-bimodule;
\medskip

\INIT $B:= \{{\pmb{f_i}} \otimes {\pmb{1}}\}_{i=1}^t \cup \{{\pmb{x^{(\epsilon_j,k)}}} \otimes {\pmb{1}}-
{\pmb{1}} \otimes {\pmb{x^{(\epsilon_j,k)}}}\}_{1\leq j \leq n,\ 1 \leq k \leq s}$;
\vskip0.1cm
\STATE  Using the Left Syzygy Module Algorithm, compute  in the PBW algebra $(\env)^{s}$
a generator system $H=\{{\pmb{h_1}},...,{\pmb{h_r}}\}$ of $Syz^l(B)$ as a left $\env$-module;
\medskip
\STATE If ${\pmb{h_i}}=({{\pmb{h_i}}}',{{{\pmb{h_i}}}}'')$ where
${\pmb{h_i}}'\in (\env)^t$ and ${{\pmb{h_i}}}''\in (\env)^{sn}$ for $1 \leq i \leq r$, take

$H:=\{{\pmb{h_1}}',...,\: {\pmb{h_r}}'\}$.
\medskip
\end{algorithmic}
\end{alg}
\begin{exmp}\label{ejemploBimoduloSicigia}
Let $R$ be the quantum plane with the PBW algebra
structure $\mathbb{C}\{x,y;\ \{yx-qxy\}, \preceq_{(2,1)}\}$,
(where $q=i\sqrt{2}$), and consider the order POT in $R^2$.
Let $F=\{(x+1,y),(xy,0)\} \subset R^2$.

Algorithm \ref{algoritmoBimoduloSicigia} takes 23.622 seconds to compute
the $R$-bimodule generator system $H$ of
$Syz(F)$ consisting of 8 elements:
\begin{eqnarray*}
\lefteqn{H = \{\: (1\otimes y-y \otimes 1, (-1+i\sqrt{2})1 \otimes 1),} \\
& & ((\frac{i\sqrt{2}}{6}-\frac{1}{3})y \otimes x +(-\frac{2}{3}-\frac{i\sqrt{2}}{3})xy \otimes 1,
1\otimes x+1 \otimes1),\: (0,1\otimes y+\frac{i\sqrt{2}}{2}y \otimes 1),\\
& & ((\frac{i\sqrt{2}}{6}-\frac{1}{3})y \otimes x +(-\frac{2}{3}+\frac{i\sqrt{2}}{3})xy \otimes 1, i\sqrt{2}x\otimes 1+1 \otimes1),\\
& & ((-\frac{1}{3}-\frac{i\sqrt{2}}{3})y \otimes x +(-\frac{2i\sqrt{2}}{3}-\frac{2}{3})xy \otimes 1, i\sqrt{2}\otimes x+i\sqrt{2} \otimes1),\\
& & (-y \otimes y+y^2\otimes 1+ i\sqrt{2}\otimes y, -i\sqrt{2}y \otimes 1),\\
& & ((\frac{1}{3}+\frac{i\sqrt{2}}{3})1 \otimes x^2 +(\frac{1}{3}-\frac{2i\sqrt{2}}{3})x \otimes x+(-\frac{2}{3}+\frac{i\sqrt{2}}{3})x^2 \otimes 1, 0),\\
& & ((\frac{1}{3}+\frac{i\sqrt{2}}{3})1 \otimes xy +(\frac{i\sqrt{2}}{6}-\frac{1}{3})y \otimes x+(-\frac{1}{3}-\frac{i\sqrt{2}}{3})x \otimes y+
(\frac{1}{3}+\frac{i\sqrt{2}}{3})xy \otimes 1, 0)\: \}.
\end{eqnarray*}
\end{exmp}
Although elimination techniques are useful at solving several problems in
Mo\-du\-le Theory,
they are computationally inefficient, mainly because elimination
orders are unavoidably used.

On the other hand, it has been proved, first in the commutative case
and then using left syzygy $R$-modules
being $R$ a non-commutative ring (such as a PBW algebra),
that syzygies provide a much more efficient treatment,
for example, in the computation of the intersection of left $R$-submodules of $R^s$,
ideal quotients,
kernels of homomorphisms of left $R$-submodules, etc.\ \citep[see][]{BGV}.

In what follows, we will see that some applications of left syzygies
can be generalized using the new definition of syzygy bimodules, so that, for example, it
is possible to give an algorithm to compute a finite intersection of $R$-subbimodules
of $R^s$ when, as natural, two-sided input data are given. Further applications will be
studied in the sequels.

The following result states a general property which the above-mentioned algorithm
will be based on.

\begin{lem}\label{propPrincipal}
Let $M$ be an $R$-subbimodule of $R^s$ such that there exist $p, q \geq 1$ and
$H=[ \frac{H_1}{H_2}] \in M_{(s+p)\times q} (R)$, where $H_1 \in M_{s\times q}(R)$ and
$H_2 \in M_{p\times q}(R)$, satisfying the following two conditions:

i) \qquad $(\mathfrak{m}^s({\pmb{h}}) \otimes {\pmb{1}})H_1={\pmb{h}}H_1,\ \forall {\pmb{h}} \in (\env)^s$;

ii) \qquad $M=\{\pmb{{\pmb{h}}} \in R^s\: ;\: \exists {{\pmb{h}}}'' \in (\env)^p\ {\rm{such\ that\ }}
({\pmb{h}} \otimes {\pmb{1}}, {\pmb{h}}'') \in Syz(H)\}$.

Let us split up each element ${\pmb{h}} \in (\env)^{s+p}$ into
${\pmb{h}}=({\pmb{h}}',{{\pmb{h}}}'')$ with ${{\pmb{h}}}' \in (\env)^s$ and
${{\pmb{h}}}'' \in (\env)^p$.
\begin{enumerate}
\item If $\{{\pmb{h_1}},...,\: {\pmb{h_t}}\} \subseteq (\env)^{s+p}$ is an
$R$-bimodule generator system of $Syz(H)$ then
$M=\: _R\langle \mathfrak{m}^s({{\pmb{h_1}}}'),...,\mathfrak{m}^s({{\pmb{h_t}}}')\rangle_R$.

\item Furthermore, if $\{{\pmb{h_1}},...,\: {\pmb{h_t}}\}$ is a left
Gr\"obner basis of $Syz(H)$ (as a left $\env$-module) for the order POT
in $(\env)^{s+p}$
and
any of 
$\preceq^{\ast},\preceq^c,\preceq_{\ast},\preceq_c$ in $\env$, then
$\{\mathfrak{m}^s({{\pmb{h_1}}}'),...,\mathfrak{m}^s({{\pmb{h_t}}}')\}\setminus \{0\}$
is a two-sided Gr\"obner basis of $M$ for POT.
\end{enumerate}
\end{lem}

\begin{thm}\label{lemaInterseccion}
Let $\{M_i\}_{i=1}^r$ be a family of $R$-subbimodules of $R^s$ and suppose that
$M_i=\: _R\langle {\pmb{f_1^i}},...,\: {\pmb{f_{t_i}^i}} \rangle_R \subseteq R^s$. Then
\[
\bigcap_{i=1}^{r} M_i=\{{\pmb{h}} \in R^s\: ; \: \exists {{\pmb{h}}}'' \in (\env)^{\sum_{j=1}^r{t_j}}
{\rm{\ such\ that\ }} ({\pmb{h}} \otimes1, {{\pmb{h}}}'') \in Syz(H)\},
\]
where
\[
H=\left(
\begin{array}{ccc}
I_s & \cdots & I_s \\
{\pmb{f_1^1}} & \cdots & {\pmb{0}} \\
\vdots & & \vdots \\
{\pmb{f_{t_1}^1}} & \cdots & {\pmb{0}} \\
\vdots & & \vdots \\
{\pmb{0}} & \cdots & {\pmb{f_1^r}} \\
\vdots & & \vdots \\
{\pmb{0}} & \cdots & {\pmb{f_{t_r}^r}} \\
\end{array}\right) \in M_{(s+\sum_{i=1}^r t_i) \times rs}(R)
\]
\end{thm}
From \ref{propPrincipal} and \ref{lemaInterseccion}, an algorithm to
compute finite intersections of $R$-subbimodules of $R^s$ may be formulated (see Algorithm \ref{algoritmoInterseccionSubimodulos}).

\begin{cor}
Let $M_i$ and $H$ be as in \ref{lemaInterseccion}.\ If $Syz(H)=\:_R\langle {\pmb{g_1}},...,\: {\pmb{g_t}} \rangle_R$
with ${\pmb{g_k}}=({{\pmb{g_k}}}',{{\pmb{g_k}}}'') \in (\env)^{s} \times (\env)^{\sum_{j=1}^r t_j}$ for all $1 \leq k \leq t$,
then
\[
\bigcap_{i=1}^{r} M_i\ =\: _R\langle \mathfrak{m}^s({{\pmb{g_1}}}'),...,\mathfrak{m}^s({{\pmb{g_t}}}')\rangle_R.
\]
If $G=\{{\pmb{g_1}},...,\: {\pmb{g_t}}\}$
is also a left Gr\"obner basis of $Syz(H)$ (as left $\env$-module) with
POT in $(\env)^{s+\sum_{j=1}^r t_j}$,
then $\{\mathfrak{m}^s({{\pmb{g_1}}}'),...,\: \mathfrak{m}^s({{\pmb{g_t}}}')\}
\setminus \{0\}$ is a two-sided Gr\"obner basis
of $\bigcap_{i=1}^{r} M_i$ with POT in $R^s$.
\end{cor}

\begin{alg}
Intersection of $R$-subbimodules of $R^s$
\begin{algorithmic}\label{algoritmoInterseccionSubimodulos}
\smallskip
\REQUIRE $\{M_i\}_{i=1}^r$, a family of $R$-subbimodules of $R^s$ with
$M_i=\: _R\langle {\pmb{f_1^i}},...,\: {\pmb{f_{t_i}^i}} \rangle_R$;
\ENSURE $M$, a finite generator system of $\bigcap_{i=1}^{r} M_i$ as an $R$-bimodule;
\medskip

\INIT
\[
H:=\left(
\begin{array}{ccc}
I_s & \cdots & I_s \\
{\pmb{f_1^1}} & \cdots & {\pmb{0}} \\
\vdots & & \vdots \\
{\pmb{f_{t_1}^1}} & \cdots & {\pmb{0}} \\
\vdots & & \vdots \\
{\pmb{0}} & \cdots & {\pmb{f_1^r}} \\
\vdots & & \vdots \\
{\pmb{0}} & \cdots & {\pmb{f_{t_r}^r}} \\
\end{array}\right) \in M_{(s+\sum_{j=1}^r t_j)\times rs}(R);
\]
\STATE  Using the Syzygy Bimodule Algorithm compute a
generator system $G=\{{\pmb{g_1}},...,{\pmb{g_t}}\}$
of $Syz(H)$ as an $R$-bimodule;
\smallskip
\STATE If ${\pmb{g_k}}=(\pmb{{g_k}}',{\pmb{g_k}}'')$ where ${\pmb{g_k}}' \in (\env)^s$ and
${\pmb{g_k}}'' \in (\env)^{\sum_{j=1}^r t_j}$ for $1 \leq k \leq t$, \\

\smallskip
take $M:=\{\mathfrak{m}^s({\pmb{g_1}}'),...,\mathfrak{m}^s({\pmb{g_t}}')\}$.
\medskip

\end{algorithmic}
\end{alg}
%
%
\begin{exmp}
Let $R$ be the quantum plane (as in example~\ref{ejemploBimoduloSicigia})
and consider the order POT in $R^2$.
Let $M_1$ and $M_2$ be the $R$-subbimodules of $R^2$ generated by
$\{(2x^2+2x,-y),(0,-8),(-3xy,0)\}$ and
$\{(x+2,0),(1,-y)\}$, respectively.

Algorithm \ref{algoritmoInterseccionSubimodulos} takes 111.418 seconds
to compute the $R$-bimodule generator system $M$ of
$M_1 \cap M_2$, consisting of 8 elements:
\begin{eqnarray*}
\lefteqn{M = \{\: (\frac{4i\sqrt{2}}{3}x^2y+\frac{7}{3}xy,\frac{i\sqrt{2}}{3}y^2),\:
(-\frac{2i\sqrt{2}}{3}x^3+(-\frac{2i\sqrt{2}}{3}-\frac{4}{3})x^2-\frac{4}{3}x,\frac{i\sqrt{2}}{3}xy+\frac{2}{3}y), } \\
& & (-\frac{5}{3}x^2y+\frac{4i\sqrt{2}}{3}xy,-\frac{2}{3}y^2),\:
((-1-\frac{5i\sqrt{2}}{3})x^3-(\frac{19}{3}+\frac{5i\sqrt{2}}{3})x^2-\frac{16}{3}x, \\
& & \frac{4i\sqrt{2}}{3}xy+\frac{8}{3}y),\:
(5x^2y+(-1-3i\sqrt{2})xy,2y^2),\: (\frac{2}{3}x^2+\frac{2}{3}x,-\frac{4}{3}y),\\
& & (\frac{5i\sqrt{2}}{3}x^2y+\frac{8}{3}xy,\frac{2i\sqrt{2}}{3}y^2),\:
(\frac{8}{3}x^2y^2-\frac{7}{3}xy^2,-\frac{1}{3}y^3),\: (-\frac{2}{3}x^3-2x^2-\frac{4}{3}x,\\
& & \frac{1}{3}xy+\frac{2}{3}y),\:
((\frac{4i\sqrt{2}}{3}+\frac{8}{3})x^2y+(\frac{4}{3}-\frac{4i\sqrt{2}}{3})xy,(\frac{i\sqrt{2}}{3}+\frac{2}{3})y^2),\\
& & (\frac{-4i\sqrt{2}}{3}x^3y-4x^2y+\frac{4i\sqrt{2}}{3}xy,-\frac{i\sqrt{2}}{3}xy^2-\frac{2}{3}y^2),\:
((-\frac{2}{3}+\frac{2i\sqrt{2}}{3})x^3+ \\
& & (-\frac{2}{3}+\frac{2i\sqrt{2}}{3})x^2,0),\: (-2x^2y+(-\frac{2}{3}+\frac{2i\sqrt{2}}{3})xy,0),\: (\frac{2}{3}x^2+\frac{2}{3}x,-\frac{4}{3}y)\: \}.
\end{eqnarray*}
\end{exmp}

\end{section}

\begin{ack}
Authors are grateful to Prof.\ Francisco Castro for his
motivation, and to Prof.\ Alain Verschoren for his corrections
and suggestions.
\end{ack}
\bibliographystyle{ams-alpha}

\end{document}